\magnification=1200
\overfullrule=0pt
\centerline {\bf A characterization related to a two-point boundary value problem}\par
\bigskip
\bigskip
\centerline {BIAGIO RICCERI}\par
\bigskip
\bigskip
\centerline {\it Dedicated to Professor Sompong Dhompongsa on his 65th birthday}\par
\bigskip
\bigskip
{\bf Abstract:} In this short note, we establish the following result: Let 
$f:[0,+\infty[\to [0,+\infty[$, $\alpha:[0,1]\to ]0,+\infty[$ be two continuous functions, with $f(0)=0$.
Assume that, for some $a>0$, the function $\xi\to {{\int_0^{\xi}f(t)dt}\over {\xi^2}}$ is 
non-increasing in $]0,a]$.\par
Then, the following assertions are equivalent:\par
\noindent
$(i)$\hskip 5pt for each $b>0$, the function $\xi\to {{\int_0^{\xi}f(t)dt}\over {\xi^2}}$ is not constant in $]0,b]$\ ;\par
\noindent
$(ii)$\hskip 5pt for each $r>0$, there exists an open interval $I\subseteq ]0,+\infty[$ such that, for every
$\lambda\in I$, the problem 
$$\cases {-u''=\lambda\alpha(t)f(u) & in $[0,1]$\cr & \cr u>0 & in $]0,1[$\cr & \cr
u(0)=u(1)=0\cr}$$
has a solution $u$ satisfying
$$\int_0^1|u'(t)|^2dt<r\ .$$
\bigskip
\indent
{\bf Key words:} Positive solutions, two-point boundary value problem, variational methods.\par
\bigskip
{\bf Mathematics Subject Classification}: 34B09, 34B18, 47J30.\par
\bigskip
\bigskip
The aim of this very short note is to establish a characterization concerning the problem
$$\cases {-u''=\lambda\alpha(t)f(u) & in $[0,1]$\cr & \cr u>0 & in $]0,1[$\cr & \cr
u(0)=u(1)=0\cr}\eqno{(D)}$$
where $f:[0,+\infty[\to [0,+\infty[$, $\alpha:[0,1]\to ]0,+\infty[$ are continuous functions, with $f(0)=0$,
and $\lambda>0$.\par
\smallskip
For each $\xi\geq 0$, set
$$F(\xi)=\int_0^{\xi}f(t)dt\ .$$
Here is our result:\par
\medskip
THEOREM 1. - {\it Assume that, for some $a>0$, the function $\xi\to {{F(\xi)}\over {\xi^2}}$ is 
non-increasing in $]0,a]$.\par
Then, the following assertions are equivalent:\par
\noindent
$(i)$\hskip 5pt for each $b>0$, the function $\xi\to {{F(\xi)}\over {\xi^2}}$ is not constant in $]0,b]$\ ;\par
\noindent
$(ii)$\hskip 5pt for each $r>0$, there exists an open interval $I\subseteq ]0,+\infty[$ such that, for every
$\lambda\in I$, problem $(D)$ has a solution $u$ satisfying
$$\int_0^1|u'(t)|^2dt<r\ .$$}
\medskip
Let $(X,\langle\cdot,\cdot\rangle)$ be a real Hilbert space. 
For each $r>0$, set
$$B_{r}=\{x\in X : \|x\|^2\leq r\}\ .$$
The key tool in our proof of Theorem 1 is provided by the following result which is
entirely based on the very recent [1]:
\medskip
THEOREM A. - {\it Let $J:X\to {\bf R}$ be a sequentially
weakly upper semicontinuos and G\^ateaux differentiable functional, with $J(0)=0$. Assume that, for some $r>0$,
there exists a global maximum $\hat x$ of $J_{|B_{r}}$ such that
$$\langle J'(\hat x),\hat x\rangle<2J(\hat x)\ .$$
Then, there exists an open interval $I\subseteq ]0,+\infty[$ such that, for every $\lambda\in I$,
the equation
$$x=\lambda J'(x)$$
has a non-zero solution lying in $\hbox {\rm int}(B_{r})$.}\par
\smallskip
PROOF. Set
$$\beta_{r}=\sup_{B_{r}}J\ ,$$
$$\delta_{r}=\sup_{x\in B_{r}\setminus\{0\}}{{J(x)}\over
{\|x\|^2}}$$
and
$$\eta(s)=\sup_{y\in B_{r}}{{r-\|y\|^2}\over {s-J(y)}}$$
for all $s\in ]\beta_{r},+\infty[$.
From Proposition 2 of [1], it follows that
$${{\beta_{r}}\over {r}}<\delta_{r}\ .$$
As a consequence, by Theorem 1 of [1], for each $s\in ]\beta_{r},r\delta_{r}[$,
the equation 
$$x={{\eta(s)}\over {2}}J'(x)$$
has a non-zero solution lying in int$(B_{r})$.  From Theorem 1 of [1] again,
we know that the function $\eta$  is convex and decreasing in $]\beta_{r},+\infty[$.
As a  consequence, the set $\eta(]\beta_{r},r\delta_{r}[)$ is an open interval.
So, the conclusion is satisfied taking
$$I={{1}\over {2}}\eta(]\beta_{r},r\delta_{r}[)$$
and the proof is complete.\hfill $\bigtriangleup$\par
\medskip
Now, we are able to prove Theorem 1.\par
\medskip
{\it Proof of Theorem 1.} We adopt the variational point of view. So, let $X$ be
the space $H^1_0(0,1)$ with the usual inner product
$$\langle u, v\rangle=\int_0^1 u'(t)v'(t)dt\ .$$
Extend the definition of $f$ (and of $F$ as well) putting it zero in $]-\infty,0[$.
Let $J:X\to {\bf R}$ be the functional defined by setting
$$J(u)=\int_0^1\alpha(t)F(u(t))dt$$
for all $u\in X$. By classical results, $J$ is $C^1$ and sequentially
weakly continuous, and (since $f\geq 0$) the solutions of problem $(D)$ are exactly the non-zero 
solutions in $X$ of the equation
$$u=\lambda J'(u)\ .$$
Let us prove that $(i)\to (ii)$. First of all, observe that, since $\xi\to {{F(\xi)}\over {\xi^2}}$ is 
non-increasing in $]0,a]$, we have
$$f(\xi)\xi\leq 2F(\xi)\eqno{(1)}$$
for all $\xi\in ]0,a]$.
Now, fix $r\in ]0,a^2]$. Since
$$\sup_{u\in X}{{\max_{[0,1]}|u|}\over {\|u\|}}\leq {{1}\over {2}}\ ,\eqno{(3)}$$
from $(1)$ it follows that
$$f(u(t))u(t)\leq 2F(u(t))\eqno{(2)}$$
for all $u\in B_{r}$ and for all $t\in [0,1]$. Now, let $u\in B_{r}$, with $\sup_{[0,1]}u>0$.
Observe that
$$\{t\in [0,1]: f(u(t))u(t)<2F(u(t))\}\neq\emptyset\ .\eqno{(4)}$$
Indeed, otherwise, in view of $(2)$ we would have
$$f(u(t))u(t)=2F(u(t))$$
for all $t\in [0,1]$ and so the function $\xi\to {{F(\xi)}\over {\xi^2}}$ would be constant
in the interval $]0,\sup_{[0,1]}u]$, against $(i)$. Then, since
 $\alpha$ is positive in $[0,1]$, from $(4)$ we infer that
$$\int_0^1\alpha(t)f(u(t))u(t)dt<2\int_0^1\alpha(t)F(u(t))dt\ .$$
This inequality can be rewritten as
$$\langle J'(u),u\rangle<2J(u)\ .$$
Therefore, all the assumptions of Theorem A are satisfied and $(ii)$ follows directly
from it.\par
Now, let us prove that $(ii)\to (i)$. Arguing by contradiction, assume that
there are $b,c>0$ such that 
$$F(\xi)=c\xi^2$$
and hence
$$f(\xi)=2c\xi$$
for all $\xi\in [0,b]$. Fix $r\in ]0,b^2]$. By $(ii)$, there exists an open interval
$I$ such that, for every $\lambda\in I$, problem $(D)$ has a solution $u$
satisfying
$$\int_0^1|u'(t)|^2dt<r\ .$$
In view of $(3)$, we have 
$$\max_{[0,1]}u\leq b$$
and so
$$f(u(t))=cu(t)$$
for all $t\in [0,1]$. In other words, for every $\lambda\in I$,
the problem
$$\cases {-u''=\lambda c\alpha(t)u & in $[0,1]$\cr & \cr u>0 & in $]0,1[$\cr & \cr
u(0)=u(1)=0\cr}$$
would have a solution. This contradicts the classical fact that the above problem has
a solution only for countably many $\lambda>0$.\hfill $\bigtriangleup$
\medskip
REMARK 1. - It is was worth noticing the following wide class of functions $f$
for which Theorem 1 applies. Namely, assume that $f$ is $2k+1$ times derivable (in a right neighbourhood
of $0$)
and that $f^{(2k)}(0)<0$ and $f^{(2m)}(0)=0$ for all $m=1,...,k-1$ if $k\geq 2$.
Then, there exists some $a>0$ such that the function $\xi\to {{F(\xi)}\over {\xi^2}}$
is decreasing in $]0,a]$. Indeed, if we put
$$\varphi(\xi)=2F(\xi)-\xi f(\xi)\ ,$$
we have $\varphi^{(2m)}(\xi)=-\xi f^{(2m)}(\xi)$ and $\varphi^{(2m+1)}(\xi)=
-f^{(2m)}(\xi)-\xi f^{(2m+1)}(\xi)$ for all $m=1,...,k$. Hence,
$\varphi(0)=\varphi^{(m)}(0)=0$ for all $m=1,...,2k$ and $\varphi^{(2k+1)}(0)>0$.
This clearly implies that, for some $a>0$, one has $\varphi(\xi)>0$ for all 
$\xi\in ]0,a]$, as claimed.

\vfill\eject
\centerline {\bf References}\par
\bigskip
\bigskip
\noindent 
[1]\hskip 5pt B. RICCERI, {\it A note on spherical maxima sharing the
same Lagrange multiplier}, preprint.\par
\bigskip
\bigskip
\bigskip
\bigskip
Department of Mathematics\par
University of Catania\par
Viale A. Doria 6\par
95125 Catania, Italy\par
{\it e-mail address}: ricceri@dmi.unict.it

\bye